\documentclass[11pt]{article}
\title
{On the number of Hamilton cycles in pseudo-random graphs}
\author{Michael Krivelevich
\thanks{School of Mathematical Sciences, Raymond and Beverly
Sackler Faculty of Exact Sciences, Tel Aviv University, Tel Aviv,
69978, Israel. Email: krivelev@post.tau.ac.il. Research supported in
part by a USA-Israel BSF grant and by a grant from the Israel
Science Foundation.}
 }

\usepackage{amsmath, amsfonts}

\begin{document}
\bibliographystyle{plain}
\maketitle
\newtheorem{thm}{Theorem}%[section]
\newtheorem{defin}{Definition}
\newtheorem{thmtool}{Theorem}[section]
\newtheorem{corollary}[thmtool]{Corollary}
\newtheorem{lem}[thmtool]{Lemma}
\newtheorem{prop}[thmtool]{Proposition}
\newtheorem{clm}[thmtool]{Claim}
\newtheorem{conjecture}{Conjecture}
\newtheorem{problem}{Problem}
\newcommand{\Proof}{\noindent{\bf Proof.}\ \ }
\newcommand{\Remarks}{\noindent{\bf Remarks:}\ \ }
\newcommand{\Remark}{\noindent{\bf Remark:}\ \ }
%Probability
\newcommand{\whp}{{\bf whp}}
\newcommand{\prob}{probability}
\newcommand{\rn}{random}
\newcommand{\rv}{random variable}
%Hypergraphs
\newcommand{\hpg}{hypergraph}
\newcommand{\hpgs}{hypergraphs}
\newcommand{\subhpg}{subhypergraph}
\newcommand{\subhpgs}{subhypergraphs}
%Letters
\newcommand{\bH}{{\bf H}}
\newcommand{\cH}{{\cal H}}
\newcommand{\cT}{{\cal T}}
\newcommand{\cF}{{\cal F}}
\newcommand{\cD}{{\cal D}}
\newcommand{\cC}{{\cal C}}

\begin{abstract}
We prove that if $G$ is an $(n,d,\lambda)$-graph (a $d$-regular
graph on $n$ vertices, all of whose non-trivial eigenvalues are at
most $\lambda)$ and the following conditions are satisfied:
\begin{enumerate}
\item $\frac{d}{\lambda}\ge (\log n)^{1+\epsilon}$ for some constant
 $\epsilon>0$;
\item $\log d\cdot \log\frac{d}{\lambda}\gg \log n$,
\end{enumerate}
then the number of Hamilton cycles in $G$ is
$n!\left(\frac{d}{n}\right)^n(1+o(1))^n$.
\end{abstract}

\section{Introduction}

The goal of this paper is to estimate the number of Hamilton cycles
in pseudo-random graphs. Putting it informally, we prove here that a
pseudo-random graph contains the right asymptotically number of
Hamilton cycles, when scaled appropriately.

Of course, the above sentence is not quite a mathematical statement,
and several of its ingredients should be explained and formalized.
The goal of this section is to provide a formal footing for this
claim.

First of all, what is the right (asymptotically) number of Hamilton
cycles? We will have to define yet the formal notion of a
pseudo-random graph to work with, but intuitively a pseudo-random
graph $G$ on $n$ vertices with $m$ edges should be similar, in some
well defined quantitative aspects, to a truly random graph on the
same number of vertices with the same (expected) number of edges. If
so, the right benchmark for the number of Hamilton cycles should
come from the standard models of random graphs.

There are quite a few available models of random graphs, of which
the most widely studied and relevant to our subject are the models
$G(n,p)$, $G(n,m)$ and $G_{n,d}$. Since over the years the random
graphs have become a part of the standard combinatorial lexicon, we
will be rather brief in defining these models. The model $G(n,p)$ of
binomial random graphs is obtained by taking $n$ labeled vertices
$1,\ldots,n=[n]$, and for each pair $1\le i<j\le n$, making $(i,j)$
into an edge independently and with probability $p=p(n)$. The
probability space $G(n,m)$ is composed of all graphs $G$ with vertex
set $[n]$ and exactly $m$ edges, where all such graphs are
equiprobable: $Pr[G]=\left(\binom{\binom{n}{2}}{m}\right)^{-1}$. The
probability space $G_{n,d}$ (assuming that the product $dn$ is even)
is composed of all $d$-regular graphs with vertex set $[n]$ and the
uniform probability measure. As customarily, we will use these
notations to denote both the corresponding probability space and a
random graph drawn from it. The random graphs $G(n,p)$ and $G(n,m)$
are quite similar to each other under proper parametrization, which
is to set $m=\binom{n}{2}p$, formal statements are available to
quantify this similarity. The random $d$-regular graph $G_{n,d}$ is
expected to resemble the binomial random graph $G(n,d/n)$ for large
enough $d=d(n)$, with a variety of concrete statements obtained to
support this paradigm. We will not dwell anymore on these concepts
and their relationships, instead referring the reader to the
standard sources in the theory of random graphs \cite{Bol-book},
\cite{JLR}.

As mentioned above, the typical number of Hamilton cycles in random
graphs will serve as a guiding line for the current research.
Consider the probability spaces $G(n,p)$ and $G(n,m)$. The number of
Hamilton cycles in the complete graph $K_n$ on $n$ vertices is
$(n-1)!/2$. Using the linearity of expectation, we obtain
immediately that if $X$ is the random variable counting the number
of Hamilton cycles, then the expectation of $X$ in the model
$G(n,p)$ is $\frac{(n-1)!}{2}p^n$, and in the model $G(n,m)$ we
have:
$$
E[X]=\frac{(n-1)!}{2}\,\frac{\binom{\binom{n}{2}-n}{m-n}}{\binom{\binom{n}{2}}{m}}\
.
$$
The above two expressions for the expectation are asymptotically
equal under the choice $m=\binom{n}{2}p$, assuming $m$ is not too
small. As the common intuition for random graphs may suggest, we
expect the random variable to be concentrated around its mean,
perhaps after some normalization (it is easy to see that the above
expressions for the expectation become exponentially large in $n$
already for $p$ inverse linear in $n$).

The reality appears to confirm this intuition --  to a certain
extent. Janson \cite{Jan94} investigated the number of Hamilton
cycles in the probability spaces $G(n,p)$ and $G(n,m)$, here are his
findings in a nutshell. As before we denote (with some ambiguity of
notation) by $X$ the random variable counting the number of Hamilton
cycles in the corresponding probability space. In the model
$G(n,m)$, assuming that $m\gg n^{3/2}$ and $\binom{n}{2}-m\gg n$,
and denoting $p=m\left/\binom{n}{2}\right.$, one has:
\begin{eqnarray*}
E[X] &=&
\frac{(n-1)!}{2}p^n\exp\left(-\frac{1-p}{p}+O\left((1-p)\frac{n^3}{m^2}\right)\right)\,,\\
Var[X] &\sim& \frac{n^3}{2m^2}(1-p)^2E^2[X]\,.
\end{eqnarray*}
and the standartized random variable $X^*=(X-E[X])/(Var[X])^{1/2}$
converges in distribution to a standard normal distribution. In
particular, for the regime $n^{3/2}\ll m \le 0.99\binom{n}{2}$, the
random variable $X$ is indeed concentrated around its expectation.
The situation appears to change around $m=\Theta(n^{3/2})$, where
the asymptotic distribution becomes log-normal instead (meaning that
$\log X$ becomes normally distributed asymptotically).

The picture in the probability space $G(n,p)$ is more involved
according to \cite{Jan94}. There we have, assuming that
$p\rightarrow \pi<1$ and $\liminf pn^{1/2}>0$:
$$
E[X] = \frac{(n-1)!}{2}p^n\,,
$$
and
$$
p^{1/2}\left(\log X -\log E[X]+\frac{1-p}{p}\right)\xrightarrow{d}
N\left(0,2(1-\pi)\right)\,.
$$

It is instructive to observe that in the latter case, and assuming
that $p\rightarrow 0$, the distribution of $X$ is in fact
concentrated way below its expectation, in particular implying that
$X/E[X]\xrightarrow{p}0$. This can be attributed to the heavy tail
of a log-normal distribution. Notice also that the number of
Hamilton cycles is more concentrated in $G(n,m)$ compared to
$G(n,p)$; this is not surprising as $G(n,m)$ is obtained from
$G(n,p)$ by conditioning on the number of edges of $G$ being exactly
equal to $m$, resulting in reducing the variance.

Though the above conclusions for the two probability spaces
$G(n,p)$, $G(n,m)$ differ quite substantially, we can put them under
one common roof by stating that (under some conditions on the
(expected) density of the random graph) one has: the number $X$ of
Hamilton cycles in a random graph with density $p$ satisfies with
high probability:
\begin{equation}\label{for1}
X = n!p^n(1+o(1))^n\ .
\end{equation}
For random graphs of density $p\ll n^{-1/2}$ not much appears to be
known about the asymptotic behavior of the number of Hamilton cycles
in corresponding random graphs. (We can mention though the result of
Cooper and Frieze \cite{CF89}, who proved that in the random graph
process typically at the very moment the minimum degree becomes two,
not only the graph is Hamiltonian but it has $(\log n)^{(1-o(1))n}$
Hamilton cycles.) This, together with the notable difference between
the results for $G(n,m)$ and $G(n,p)$, and the somewhat unexpected
form of the result in the case of $G(n,p)$, indicate that even for
the case of random graphs the question about counting the number of
Hamilton cycles is far from being trivial.

 For the probability space of random regular graphs, it is the
opposite case of sparse graphs that is relatively well understood.
Janson \cite{Jan95}, following the previous work of Robinson and
Wormald \cite{RW92}, \cite{RW94}, described the asymptotic
distribution of the number of Hamilton cycles in a random
$d$-regular graph $G_{n,d}$ for a {\em constant} $d\ge 3$. The
expression obtained is quite complicated, and we will not reproduce
it here. No results seem to be available in the literature for the
case of growing degree $d=d(n)$.

Now that we have covered briefly results about the number of
Hamilton cycles in random graphs, we switch to the pseudo-random
case, and more generally, to deterministic results. Frieze treated
the case of dense pseudo-random graphs in \cite{Fri00}. He proved
that if $G$ is a graph on $n$ vertices, meeting the following two
conditions:
 \begin{enumerate}
\item all degrees in $G$ are in the range $[d-\epsilon n,d+\epsilon n]$;
\item for every pair of disjoint sets $S,T\subset V(G)$,
$|S|,|T|\ge\epsilon n$, the number of edges between $S$ and $T$ in
$G$ in the range
$\left[\left(\frac{d}{n}-\epsilon\right)|S||T|,\left(\frac{d}{n}+\epsilon\right)|S||T|\right]$,
\end{enumerate}
then the number of Hamilton cycles in $G$ is in the range
$\left[\left(\frac{d}{n}-2\epsilon\right)^nn!,\left(\frac{d}{n}+2\epsilon\right)^nn!\right]$,
quite in line with the general paradigm (\ref{for1}). The above
assumptions on $G$ are obviously shaped after the binomial random
graph $G(n,\frac{d}{n})$; also, they tacitly assume that the typical
degree $d=d(n)$ in $G$ is linear in $n$ and $\epsilon\ll d/n$, as
taking $d=o(n)$ and $\epsilon$ constant renders both the assumptions
and the conclusion essentially meaningless. Frieze also obtained
similar results for the numbers of perfect matchings and of spanning
trees in such a pseudo-random graph in the same paper. Recently
Cuckler and Kahn \cite{CK09-1}, \cite{CK09-2} considered the case of
dense graphs. They proved that if $G$ is a graph on $n$ vertices
with the minimum degree $\delta(G)\ge \frac{n}{2}$, then not only
$G$ is Hamiltonian as asserted by the very well known Dirac theorem,
but it contains in fact at least $n!\left(\frac{1}{2}+o(1)\right)^n$
Hamilton cycles; this estimate, which is essentially optimal due to
what has been stated already about the random graph $G(n,1/2)$,
improved a prior result by S\'ark\"ozy, Selkow and Szemer\'edi
\cite{SSS03}. Cuckler and Kahn proved also that if $G$ is a
$d$-regular graph on $n$ vertices with $d\ge \frac{n}{2}$ (and
without any further assumptions on the edge distribution in $G$),
then the number of Hamilton cycles in $G$ is
$n!\left(\frac{d}{n}+o(1)\right)^n$, again as predicted by
(\ref{for1}).

Now it is about time to reveal the model of pseudo-random graphs we
adopt for this paper. As we briefly mentioned, a definition of
pseudo-random graphs is supposed to capture quantitatively their
similarity to truly random graphs of the same (expected) density.
Note that not every definition of pseudo-random graphs is suitable
for our purposes; for example, the classical definition of jumbled
graphs due to Thomason \cite{Tho87} is no good for us as it allows
occasional vertex degrees to deviate substantially from the average
degree, including the possibility of the existence of isolated
vertices, thus making any estimate of the number of Hamilton cycles
virtually impossible.

Here we will use the notion of $(n,d,\lambda)$-graphs to model
pseudo-random graphs.
\begin{defin}\label{ndl}
A graph $G$ is an {\em $(n,d,\lambda)$-graph} if $G$ has $n$
vertices, is $d$-regular, and the second largest (in absolute value)
eigenvalue of its adjacency matrix is bounded from above by
$\lambda$.
\end{defin}

This is one of the most studied notions of pseudo-random graphs. It
is very convenient for our purposes as it ensures that all degrees
are equal and also allows for a very good grip on the edge
distribution in such a graph. We will provide more technical details
about $(n,d,\lambda)$-graphs in Section \ref{sec3-1}. The reader is
referred to surveys \cite{HLW06}, \cite{KS06} for a thorough
discussion of $(n,d,\lambda)$-graphs, their examples and properties.
Let us just mention here that an $(n,d,\lambda)$-graph $G$ with
$\lambda\ll d$ resembles quite closely a binomial random graph
$G(n,\frac{d}{n})$ or a random $d$-regular graph $G_{n,d}$ in many
quantitative aspects.

Of course, before embarking on estimating the number of Hamilton
cycles in $(n,d,\lambda)$-graphs we should make sure that Hamilton
cycles do exist in such graphs. Such a statement is available indeed
\cite{KS03}, we will state and discuss it in Section \ref{sec3-2}.

We can now formulate the main result of this paper.

\begin{thm}\label{th1}
For every $\epsilon>0$ and for sufficiently large $n$ the following
is true. Let  $G$ be an $(n,d,\lambda)$-graph, satisfying the
following conditions:
\begin{gather}
 \frac{d}{\lambda}\ge (\log n)^{1+\epsilon}\,,\label{cond1}\\
\log d\cdot \log\frac{d}{\lambda}\gg \log n\,.\label{cond2}
\end{gather}
Then the number of Hamilton cycles in $G$ is asymptotically equal to
$n!\left(\frac{d}{n}\right)^n(1+o(1))^n$.
\end{thm}

Let us discuss the above statement briefly. Condition (\ref{cond1})
appears to be rather mild and is only a notch above the best known
sufficient condition for Hamiltonicity in $(n,d,\lambda)$-graphs
provided by \cite{KS03}. As for condition (\ref{cond2}), for the
(rather typical) case $\lambda\le d^{1-\alpha}$ for some constant
$\alpha>0$, (\ref{cond2}) becomes $\log d\gg \log^{1/2}n$, which is
equivalent to $d=2^{\omega(\log^{1/2}n)}$. Thus, condition
(\ref{cond2}) admits graphs of subpolynomial degrees. Of course, the
main thrust of Theorem \ref{th1} is to ensure that under some rather
mild assumptions the number of Hamilton cycles in a graph $G$ is
what is predicted by (\ref{for1}). The error term $(1+o(1))^n$ is
extremely convenient and robust as in particular it allows to sweep
under the rug even very fast growing functions of $n$, like for
example $2^{\frac{n}{\log n}}$. Due to the standard estimates on
$n!$ (say, the Stirling formula) the estimate of Theorem \ref{th1}
reads as $\left(\frac{d}{e}\right)^n(1+o(1))^n$.

The remainder of this paper is organized as follows. The next
section introduces definitions and notation used in later sections.
In Section \ref{sec3} we describe the set of tools used in our main
proof. Theorem \ref{th1} is proven then in Section \ref{sec4}.
Section \ref{sec5}, the last section of the paper, is devoted to
concluding remarks.

\section{Definitions and notation}\label{sec2}

The number of Hamilton cycles in a graph $G$ is denoted by $h(G)$.
In this paper, we consider a single edge as a cycle too. A {\em
2-factor} in a graph $G$ is a collection of vertex disjoint cycles
covering all vertices of $G$. For a 2-factor $F$ in $G$, we denote
by $c(F)$ the number of cycles of length at least 3 in $F$. For a
graph $G$ and an integer $s$, we let $f(G,s)$ be the number of
2-factors in $G$ with exactly $s$ cycles; $f(G)$ is the total number
of 2-factors in $G$. For a graph $G$ and an integer $2\le k\le
|V(G)|$ we define
$$
\phi(G,k)= \max\{f(G[V_0]): V_0\subseteq V, |V_0|=k\}\,.
$$

The other notation we use is fairly standard. In particular, given a
graph $G=(V,E)$ and vertex subsets $U,W\subseteq V$, we denote by
$e_G(U,W)$ the number of edges of $G$ with one endpoint in $U$ and
another in $W$; by $e_G(U)$ the number of edges of $G$ spanned by
$U$ (thus, $e_G(U)=\frac{1}{2}e_G(U,U)$), and by $N_G(U)$ the
external neighborhood of $U$ in $G$; whenever the identity of the
graph $G$ is clear from the context, we will omit placing it in the
index of the above notations. All logarithms are natural.

As our result is asymptotic in nature, we routinely assume that the
underlying parameter $n$ (normally standing for the number of
vertices in a graph $G$ under consideration) is large enough for our
purposes.

\section{Tools}\label{sec3}
\subsection{$(n,d,\lambda)$-graphs and the expander mixing
lemma}\label{sec3-1}

As we have already declared our model of pseudo-random graphs is
$(n,d,\lambda)$-graphs. The most basic property of an
$(n,d,\lambda)$-graph is given by the following very well known
statement, bridging between graph eigenvalues and edge distribution
and sometimes called the Expander Mixing Lemma (see, e.g. Corollary
9.2.5 of \cite{AS} or Theorem 2.11 of \cite{KS06}). Let $G$ be an
$(n,d,\lambda)$-graph. Then for any two vertex subsets $S,T\subseteq
V(G)$
\begin{equation}\label{ndl-distr}
\left|e(S,T)-\frac{d}{n}|S|\,|T|\right|\le \lambda\sqrt{|S|\,|T|}\,.
\end{equation}
This formula shows obviously the quantitative similarity of the edge
distribution in an $(n,d,\lambda)$-graph $G$ to that of a binomial
random graph $G(n,p)$ with the edge probability $p=d/n$. Indeed, in
$G(n,d/n)$ we expect $\frac{d}{n}|S||T|$ edges between $S$ and $T$,
and estimate (\ref{ndl-distr}) shows that this is basically what
happens in an $(n,d,\lambda)$-graph, assuming the sets $S,T$ are
large enough, and the so called eigenvalue ratio $d/\lambda$ is
relatively large as well. The error term in (\ref{ndl-distr}) is
governed by $\lambda$;  the smaller $\lambda$ is, the better the
edge distribution fits the expected random pattern. Speaking in more
concrete terms, one can derive from (\ref{ndl-distr}) that small
sets in an $(n,d,\lambda)$-graph expand outside substantially:
$$
|N(X)|\ge \frac{(d-2\lambda)^2}{3\lambda^2}|X|
$$
for $|X|\le \frac{\lambda^2n}{d^2}$ (see, e.g., Proposition 2.3 of
\cite{KS03}), while there is always an edge between two large enough
sets: for every pair of disjoint sets $X,Y$ with $|X|,|Y|>
\frac{\lambda n}{d}$, one has $e(X,Y)>0$; indeed, in such a case by
(\ref{ndl-distr}): $e(X,Y)\ge \frac{d}{n}|X|\,|Y|-\lambda\sqrt{|X|\,
|Y|}>\frac{d}{n}\cdot\frac{\lambda n}{d}|X|\,|Y|-\lambda|X|\,
|Y|=0$.

\subsection{Hamiltonicity in $(n,d,\lambda)$-graphs}\label{sec3-2}
The paper \cite{KS03} provides a sufficient condition for
Hamiltonicity in $(n,d,\lambda)$-graphs in terms of the eigenvalue
ratio. It is proven in \cite{KS03} that if $n$ is large enough and
$$
\frac{d}{\lambda}\ge \frac{1000\log n\log\log\log n}{(\log\log
n)^2}\,,
$$
then an $(n,d,\lambda)$-graph $G$ is Hamiltonian. (A related result
is \cite{HKS09}, where a sufficient condition for Hamiltonicity of a
general graph $G$ is stated in terms of expansion and a
connectivity-type condition.)

The argument of \cite{KS03} utilizes the ingenious
rotation-extension technique of P\'osa \cite{Pos}, very frequently
used in papers on Hamiltonicity of random and pseudo-random graphs.
Since we will not apply it directly in this paper, we will skip its
detailed description, instead referring to it in general terms.

For the purposes of this paper, we need a certain, quite
straightforward, modification of the argument of \cite{KS03}. This
modification will allow us to control the number of rotations
performed when constructing a Hamilton cycle.
\begin{lem}\label{lem-KS}
For every $\epsilon>0$ there exist $C=C(\epsilon)>0$ and
$n_0=n_0(\epsilon)>0$ such that for every integer $n\ge n_0$ the
following is true. Let
\begin{equation}\label{lem-KS-cond}
\frac{d}{\lambda}\ge (\log n)^{1+\epsilon}\,.
\end{equation}
Let $G$ be an $(n,d,\lambda)$-graph. Then $G$ is connected. Let
further $P_0$ be a path in $G$. Then there is a path $P^*$ in $G$
connecting vertices $a$ and $b$ so that:
\begin{enumerate}
\item $V(P^*)=V(P_0)$;
\item $|E(P_0)\bigtriangleup E(P^*)|\le
 \frac{C\log n}{\log\frac{d}{\lambda}}$;
\item $(a,b)\in E(G)$, or $G$ contains an edge between $\{a,b\}$ and
$V(G)-V(P^*)$.
\end{enumerate}
\end{lem}

The main quantitative conclusion of the above lemma is its second
consequence above, allowing to bound the number of rotations needed
to get from $P_0$ to $P^*$. Once we obtain the path $P^*$ as in the
lemma, we can close it to a cycle -- which is either Hamiltonian, or
can be used to find a path longer than $P^*$ due to connectivity by
adjoining a vertex outside $V(P^*)$; the other alternative is to
directly extend $P^*$ to a longer path by appending a new vertex to
one of its endpoints $a,b$. Of course this shows that an
$(n,d,\lambda)$-graph satisfying (\ref{lem-KS-cond}) is Hamiltonian,
but our main point here is different -- we say that after
$O\left(\frac{\log n}{\log\frac{d}{\lambda}}\right)$ rotations/edge
changes starting from any given path we are at least one step closer
to Hamiltonicity.

For the sake of our alert readers we now indicate briefly how the
proof presented in \cite{KS03} can be adjusted to give Lemma
\ref{lem-KS}. The focus of our attention is Section 3 of that paper.
The arguments of Section 3.1 do not require any modification; the
definition of $\rho$ from that subsection stays the same. In Section
3.2 we take $k=2$. This implies $\alpha=\Theta(1)$. We get sets
$C_1,C_2$ of sizes $|C_1|,|C_2|=\Theta\left(\frac{n}{\rho}\right)=
\Theta\left(\frac{n}{\frac{\log n}{\log\frac{d}{\lambda}}}\right)$.
Then in Proposition 3.2 we find $C_1'\subseteq C_1$ with
$int(C_1')=\Theta\left(\frac{n}{\rho}\right)$ such that every vertex
$v\in C_1'$ has $\Omega\left(\frac{d}{\rho}\right)$ neighbors in
$int(C_1')$. This argument would use the estimate:
$\frac{d}{\lambda}\,\cdot\,\log\frac{d}{\lambda}\gg \log n$. We
argue similarly to find a set $C_2'\subseteq C_2$. In Proposition
3.5 we get to a set $T_i$ with $|T_i|\ge\frac{\lambda n}{d}$ in
$O\left(\frac{\log n}{\log\frac{d}{\lambda}}\right)$ rotations. The
reason is that in every induced subgraph $G_0$ of $G$ of minimum
degree $\Omega\left(\frac{d}{\rho}\right)$ every small set expands
itself outside by the factor of
$\Omega\left(\frac{d^2}{\rho^2\lambda^2}\right)$. Therefore we need
$$
O\left(\frac{\log\left(\frac{\lambda n}{d}\right)}
            {\log\left(\frac{d^2}{\rho^2\lambda^2}\right)}\right)=
            O\left(\frac{\log
            n}{\log\left(\frac{d}{\rho\lambda}\right)}\right)
$$
rotations. In order to get to $T_i$ in $O\left(\frac{\log
n}{\log\frac{d}{\lambda}}\right)$ rotations we need to require:
$\log\left(\frac{d}{\rho\lambda}\right)=\Omega\left(\log\frac{d}{\lambda}\right)$,
which is equivalent to:
$\frac{d}{\rho\lambda}\ge\left(\frac{d}{\lambda}\right)^{\delta}$
for some $\delta>0$. Recalling that $\rho=\Theta\left(\frac{\log
n}{\log\frac{d}{\lambda}}\right)$, we see that this condition is
satisfied if $\frac{d}{\lambda}\ge (\log n)^{1+\epsilon}$ -- which
is exactly assumption (\ref{cond1}) of Theorem \ref{th1}.

\subsection{Permanent estimates}\label{sec3-3}
There is a well known and frequently used connection between cycles
and cycle factors in graphs and matrix permanents. This connection
has been utilized in several papers on Hamiltonicity, see, e.g.,
\cite{FK05}, \cite{KKO11-1}, \cite{KKO11-2}. Permanent estimates
play a crucial role in our arguments too.

We need both upper and lower bounds for permanents. The upper bound,
conjectured by Minc and proved by Bregman \cite{Bre73}, together
with an elementary convexity argument (see, e.g., Corollary 3 at p.
64 of \cite{AS}), gives:

\begin{lem}\label{per-u}
Let $A$ by an $n\times n$ matrix of $0-1$ with $t$ ones altogether.
Then $per(A)\le \prod_{i=1}^n(r_i!)^{1/r_i}$, where $r_i$ are
integers satisfying $\sum_{i=1}^n r_i=t$ and as equal as possible.
\end{lem}

The lower bound, conjectured by van der Waerden and proved by
Egorychev \cite{Ego81} and by Falikman \cite{Fal81} is as follows:

\begin{lem}\label{per-i}
Let $A$ be an $n\times n$ doubly stochastic matrix. Then $per(A)\ge
\frac{n!}{n^n}$.
\end{lem}

\section{Proof of Theorem \ref{th1}}\label{sec4}

As we have indicated already we base our proof (both lower and upper
bounds) on a connection between 2-factors and Hamilton cycles in
graphs and permanents of graph matrices. Let $A$ be the adjacency
matrix of $G$. Then $A$ is an $n$-by-$n$ matrix of 0-1 with exactly
$d$ ones in each row and column, implying in particular that the
matrix $\frac{1}{d}A$ is doubly stochastic.

Consider the permanent of $A$. Each generalized diagonal
contributing to the permanent corresponds naturally to a 2-factor
(obtained by taking the edges corresponding to the entries of this
generalized diagonal); moreover, each 2-factor $F$ is counted
exactly $2^{c(F)}$ times (as there are two ways to orient each of
$c(F)$ non-trivial cycles from $F$). We thus get:
 \begin{equation}\label{per-fac}
 per(A) = \sum_{\mbox{$F$ - 2-factor in $G$}}2^{c(F)}\,.
 \end{equation}

Now, the upper bound of Theorem \ref{th1} follows immediately from
the above estimate and Bregman's theorem (Lemma \ref{per-u}):
$$
h(G)\le f(G)\le \sum_{\mbox{$F$ - 2-factor in
$G$}}2^{c(F)}=per(A)\le (d!)^{\frac{n}{d}}\,.
$$
Plugging in the estimate $d!\le d(d/e)^d$, we get
$$
h(G)\le \left(\frac{d}{e}\right)^n\,\cdot\, d^{\frac{n}{d}}=
\left(\frac{d}{e}\right)^n(1+o(1))^n\,,
$$
proving the upper bound. (Observe that the proof shows that the
obtained upper bound is valid in fact for {\em any} $d$-regular
graph $G$.)

The lower bound is much more challenging. Before delving into the
details and calculations of the proof, we outline the main steps of
our argument.
\begin{enumerate}
\item We first use (\ref{per-fac}) and the Egorychev-Falikman
theorem to estimate from below the number of 2-factors in $G$,
weighted by their numbers of cycles.
\item Then we prove that the contribution of 2-factors with many cycles
to this number is rather insignificant; here we will use estimate
(\ref{ndl-distr}) on the edge distribution of $(n,d,\lambda)$-graphs
and Bregman's theorem.
\item Then we prove that each 2-factor with relatively few cycles
can be converted into a Hamilton cycle using relatively few
rotations; here Lemma \ref{lem-KS} is applied.
\item We conclude that since we have
$\left(\frac{d}{e}\right)^n(1-o(1))^n$ 2-factors with relatively few
cycles, each being relatively close to a Hamilton cycle, the number
of Hamilton cycles should be large as well, bringing us to the
desired bound.
\end{enumerate}

We now start filling in the details of the proof. From
(\ref{per-fac}) we get:
$$
per(A) = \sum_{\mbox{$F$ - 2-factor in $G$}}2^{c(F)}\le
\sum_{s=1}^{n/2}f(G,s)\cdot 2^s\,.
$$
Applying the van der Waerden Conjecture (Lemma \ref{per-i}) to the
doubly stochastic matrix $\frac{1}{d}A$, we obtain:
\begin{equation}\label{l9}
\sum_{s=1}^{n/2}f(G,s)\cdot 2^s\ge per(A)\ge
n!\left(\frac{d}{n}\right)^n\ge \left(\frac{d}{e}\right)^n\,.
\end{equation}
Set
$$
s^*=\frac{20n}{\log^2d}\,.
$$
We will show that the contribution of 2-factors with many cycles to
the last sum is negligible:
\begin{equation}\label{l10}
\sum_{s>s^*}f(G,s)\cdot 2^s =
o\left(\left(\frac{d}{e}\right)^n\right)\,.
\end{equation}

Let $s>s^*$. Our goal is to estimate the term $f(G,s)\cdot 2^s$ from
above. Define
$$
s_1 =\frac{4s}{\log d}\,.
$$
If a 2-factor $F$ has s cycles then (by taking its shortest cycles)
we see that $F$ has $s_1$ cycles of total length $t\le
\frac{s_1}{s}n=\frac{4n}{\log d}$.

Fix $t\le \frac{4n}{\log d}$. If $(k_1,\ldots,k_{s_1})$ is a vector
of cycle lengths satisfying $\sum_{i=1}^{s_1}k_i=t$, the number of
2-factors whose $s_1$ shortest cycles are of lengths
$(k_1,\ldots,k_{s_1})$ is at most:
\begin{equation}\label{l11}
\binom{n}{s_1}\,\cdot\,\prod_{i=1}^{s_1}\frac{d^{k_i-1}}{k_i}
\,\cdot\, \phi(G,n-t) \le \binom{n}{s_1}d^{t-s_1}\phi(G,n-t)
\end{equation}
(for the expression in the left hand side above, first choose one
vertex from each cycle, then for each of the cycles construct a path
of length $k_i-1$ from the corresponding chosen vertex; once the
$s_1$ cycles are laid out, complete their union to a 2-factor
spanned by the remaining $n-t$ vertices).

Now we estimate $\phi(G,n-t)$. Let $V_0$ be a subset of $V(G)$ of
cardinality $|V_0|=t$. Denote by $A_1$ the adjacency matrix of the
subgraph $G[V-V_0]$. As we argued before, the number of 2-factors in
this graph is at most $per(A_1)$. In order to estimate $per(A_1)$,
notice that
$$
e_G(V_0)\le \frac{t^2}{2}\,\frac{d}{n}+\lambda t
$$
by estimate (\ref{ndl-distr}). It thus follows that
$$
e_G(V_0,V-V_0)=dt-2e_G(V_0)\ge dt-\frac{dt^2}{n}-2\lambda t\,.
$$
We derive:
$$
2e_G(V-V_0)=d(n-t)-e(V_0,V-V_0)\le d(n-t)-dt+\frac{dt^2}{n}+2\lambda
t\,.
$$
It thus follows that the average degree in the induced subgraph
$G[V-V_0]$ is
\begin{eqnarray*}
\frac{2e_G(V-V_0)}{n-t}&\le&
d-\frac{dt}{n-t}+\frac{dt^2}{n(n-t)}+\frac{2\lambda t}{n-t}\\
&=& d\left(1-\frac{t}{n}\right)+\frac{2\lambda t}{n-t} =: d_1\,.
\end{eqnarray*}
Then by Lemma \ref{per-u}
\begin{eqnarray*}
per(A_1)&\le& (\lceil d_1\rceil!)^{\frac{n-t}{\lfloor d_1\rfloor}}
\le \left(\left(\frac{d_1}{e}\right)^{d_1}\cdot d_1\right)
^{\frac{n-t}{\lfloor d_1\rfloor}}
 \le \left(\frac{d_1}{e}\right)^{n-t}\cdot
 \left(\frac{d_1}{e}\right)^{n\left(\frac{d_1}{\lfloor
 d_1\rfloor}-1\right)}\cdot d_1^{\frac{n}{\lfloor d_1\rfloor}}\\
 &\le& \left(\frac{d_1}{e}\right)^{n-t}\cdot d_1^{\frac{2n}{d_1}}\cdot d_1^{\frac{2n}{d_1}}
 \le\left(\frac{d_1}{e}\right)^{n-t}\cdot e^{\frac{5n\log d}{d}}\,.
\end{eqnarray*}

Substituting the expression for $d_1$ in the estimate above we get:
\begin{eqnarray*}
per(A_1)&\le& \left(\frac{d\left(1-\frac{t}{n}\right)+\frac{2\lambda
t}{n-t}}{e}\right)^{n-t}\cdot e^{\frac{5n\log d}{d}}\\
&\le& \frac{d^{n-t}}{e^{n-t}}\cdot
\left(1-\frac{t}{n}\right)^{n-t}\cdot \left(1+\frac{4\lambda
t}{d(n-t)}\right)^{n-t}\cdot e^{\frac{5n\log d}{d}}\\
&\le&
\frac{d^{n-t}}{e^n}\cdot\exp\left\{\frac{t^2}{n}+\frac{4\lambda
t}{d}+\frac{5n\log d}{d}\right\}\,.
\end{eqnarray*}

The above is an upper bound on $\phi(G,n-t)$. Plugging it into
(\ref{l11}) and estimating the number of solutions of
$k_1+\ldots+k_{s_1}=t$ in positive integers by $\binom{t+s_1}{s_1}$,
we have:
\begin{eqnarray*}
f(G,s)\cdot 2^s &\le& \sum_{t\le\frac{4n}{\log
d}}\binom{n}{s_1}\binom{t+s_1}{s_1}\left(\frac{d}{e}\right)^n\cdot
2^s\cdot\exp\left\{\frac{t^2}{n}+\frac{4\lambda t}{d}+\frac{5n\log
d}{d}\right\}\cdot d^{-s_1}\\
&\le& \left(\frac{d}{e}\right)^n\,\sum_{t\le\frac{4n}{\log d}}
\left(\frac{en}{s_1}\right)^{s_1}\,\left(\frac{5t}{s_1}\right)^{s_1}\cdot
2^s\cdot\exp\left\{\frac{t^2}{n}+\frac{4\lambda t}{d}+\frac{5n\log
d}{d}\right\}\cdot d^{-s_1}\,.
\end{eqnarray*}

The $t$-th summand in the sum above is at most
$$
\left(\frac{15nt}{ds_1^2}\right)^{s_1}\cdot 2^{\frac{s_1\log
d}{4}}\cdot\exp\left\{\frac{t^2}{n}+\frac{4\lambda
t}{d}+\frac{5n\log d}{d}\right\}
 \le \left(\frac{nt}{d^{3/4}s_1^2}\right)^{s_1}
 \exp\left\{\frac{t^2}{n}+\frac{4\lambda t}{d}+\frac{5n\log d}{d}\right\}
 \,.
$$
Since $s_1=\frac{4s}{\log d}\ge \frac{4s^*}{\log
d}=\frac{80n}{\log^3d}$, we have:
$\left(\frac{d^{3/4}s_1^2}{nt}\right)^{s_1}\ge d^{0.7s_1}\ge
e^{\frac{50n}{\log^2d}}$. For the terms in the exponent
$\exp\{\ldots\}$ above, we have the following estimates:
\begin{gather*}
t\le \frac{4n}{\log d}\Rightarrow \frac{t^2}{n}\le
\frac{16n}{\log^2d}\,,\\
\frac{\lambda t}{d}\le \frac{4n}{\log d}\cdot \frac{1}{(\log
n)^{1+\epsilon}}=o\left(\frac{n}{\log^2d}\right)\,,\\
\frac{5n\log d}{d}=o\left(\frac{n}{\log^2d}\right)\,,
\end{gather*}
and thus $\exp\{\ldots\}\le e^{\frac{17n}{\log^2d}}$. It follows
that
$$
f(G,s)\cdot 2^s\le
\left(\frac{d}{e}\right)^n\sum_{t}e^{-\frac{50n}{\log^2d}+\frac{17n}{\log^2d}}
= \left(\frac{d}{e}\right)^n\cdot o\left(\frac{1}{n}\right)\,.
$$
Hence $\sum_{s>s^*}f(G,s)\cdot 2^s=o\left((d/e)^n\right)$,
establishing (\ref{l10}). We obtain from (\ref{l9}):
\begin{equation}\label{l12}
\sum_{s\le s^*}f(G,s)\ge
\frac{1-o(1)}{2^{s^*}}\left(\frac{d}{e}\right)^n\,.
\end{equation}

Let now $F$ be a 2-factor in $G$ with $s\le s^*$ cycles. We can turn
in into a Hamilton cycle in $G$ be deleting and inserting some (few)
edges as follows. Let $C$ be an arbitrary cycle in $G$. By
connectivity one of the vertices of $C$, say, $v$ has a neighbor
outside $C$ -- unless of course $C$ is already Hamiltonian. Open $C$
up be deleting an edge of $C$ incident to $v$ (no need to do so if
$C$ is just an edge), we get a path $P$. Since there is an edge
$e\in E(G)$ between an endpoint of $P$ and some other cycle $C'$ in
$F$ we append this edge to $P$, go through it to $C'$, open it up be
deleting an edge of $C'$ incident to $e$ to get a longer path $P'$
and repeat the argument. If at some point there are no edges between
the endpoints of the current path $P''$ and other cycles from $F$,
then we can rotate $P''$ using Lemma \ref{lem-KS} to close it to a
cycle or to extend it outside. In all cases according to Lemma
\ref{lem-KS} we invest  $O\left(\frac{\log
n}{\log\frac{d}{\lambda}}\right)$ edge replacements to reduce the
number of cycles by at least 1, and thus after
$O\left(s\cdot\frac{\log n}{\log\frac{d}{\lambda}}\right)$
replacements we get to a Hamilton cycle.

Looking at it from the other side, observe that a given Hamilton
cycle $H$ in $G$ is at distance at most $k$ from at most
$\binom{n}{k}\,d^{2k}$ 2-factors in $G$ (first choose $k$ edges of
$H$ to be deleted, thus obtaining a collection of at most $k$ paths;
these paths should be then tailored into a 2-factor, and the number
of choices here is at most $d$ per each of the at most $2k$
endpoints of the paths). Hence
$$
\sum_{s\le s^*}f(G,s)\le h(G)\,\binom{n}{k}d^{2k}
$$
for $k=O\left(s^*\cdot\frac{\log n}{\log\frac{d}{\lambda}}\right)$.
We obtain from (\ref{l12}):
$$
h(G)\ge \frac{\sum_{s\le s^*}f(G,s)}{\binom{n}{k}d^{2k}}\ge
\left(\frac{d}{e}\right)^n\cdot
\frac{1-o(1)}{2^{s^*}\binom{n}{k}d^{2k}}\,.
$$
Since $k=O\left(s^*\cdot\frac{\log n}{\log\frac{d}{\lambda}}\right)
=O\left(\frac{n}{\log^2d}\cdot
 \frac{\log n}{\log\frac{d}{\lambda}}\right)=o(n)$,
 we have $\binom{n}{k}=2^{o(n)}$. Also,
$$
d^{2k}\le d^{Cs^*\,\frac{\log n}{\log\frac{d}{\lambda}}}=
e^{\frac{Cn}{\log d}\,\frac{20\log
n}{\log\frac{d}{\lambda}}}=e^{o(n)}\,,
$$
where the last estimate is due to our assumption (\ref{cond2}). It
thus follows that
$$
h(G)\ge \left(\frac{d}{e}\right)e^{-o(n)}=
n!\left(\frac{d}{n}\right)^n(1-o(1))^n\,,
$$
completing the proof of the lower bound of Theorem \ref{th1}.

\section{Concluding remarks}\label{sec5}

We have proven that an $(n,d,\lambda)$-graph $G$, a quite popular
model of pseudo-random graphs, contains $n!(d/n)^n(1+o(1))^n$
Hamilton cycles, as to be expected based on the intuition borrowed
from random graphs; this is under additional assumptions
(\ref{cond1}) and (\ref{cond2}) on the degree $d$ and the spectral
ratio $d/\lambda$. It would be nice to relax the second assumption
to make the result applicable to $d$-regular graphs on $n$ vertices
with the degree $d=d(n)$ as low as polylogarithmic in $n$. Another
attractive avenue to explore is to try and obtain similar estimates
for other models of pseudo-random graphs, perhaps less
rigid/restrictive than the model of $(n,d,\lambda)$-graphs.

Finally, let us note that our bound on the number of Hamilton cycles
in an $(n,d,\lambda)$-graph can be used to bound the number of
perfect matchings (this connection has been exploited in, e.g.,
\cite{Fri00}, \cite{CK09-2}). Let $m(G)$ denote the number of
perfect matchings in $G$. Let now $G$ be an $(n,d,\lambda)$-graph
with $n$ even. Observe that each Hamilton cycle in $G$ is a union of
two perfect matchings. This implies $h(G)\le\binom{m(G)}{2}$, and
thus a lower bound on $h(G)$ supplied a lower bound on $m(G)$. For
the upper bound, we can use for example a result of Alon and
Friedland \cite{AF08}, who proved in particular that the number of
perfect matchings in any $d$-regular graph $G$ on $n$ vertices is at
most $(d!)^{\frac{n}{2d}}$. These two bounds combined together show
that the number of perfect matchings in an $(n,d,\lambda)$-graph $G$
satisfying the conditions of Theorem \ref{th1} (with $n$ even, of
course) is asymptotic to
$\left(\frac{d}{e}\right)^{\frac{n}{2}}(1+o(1))^n$.

% \noindent{\bf Acknowledgement.}

\end{document}